\documentclass[11pt]{article}
\usepackage{CJK}
\usepackage{amsmath}
\usepackage{amsmath,amssymb,latexsym,color}
\usepackage[mathscr]{eucal}
\usepackage{graphicx}
\usepackage{pstricks}
\newtheorem{thm}{Theorem}[section]
\newtheorem{defin}[thm]{Definition}

\newtheorem{prop}[thm]{Proposition}
\newtheorem{lemma}[thm]{Lemma}
\newtheorem{cor}[thm]{Corollary}

\newtheorem{example}[thm]{Example}

\newcommand{\proof}{{\it Proof.\quad}}
\newcommand{\qed}{\hfill\Box\medskip}
\begin{document}
\begin{CJK*}{GBK}{song}

\title{\bf Identifying codes of corona product graphs}

\author{Min Feng\quad Kaishun Wang\footnote{Corresponding author. E-mail address: wangks@bnu.edu.cn}\\
{\footnotesize   \em  Sch. Math. Sci. {\rm \&} Lab. Math. Com. Sys.,
Beijing Normal University, Beijing, 100875,  China}}
\date{}
\maketitle

\begin{abstract}
For a vertex $x$ of a graph $G$, let $N_G[x]$ be the set of $x$
with all of its neighbors in $G$. A set $C$ of vertices is an {\em
identifying code} of $G$ if the sets $N_G[x]\cap C$ are nonempty and
distinct for all vertices $x$. If $G$ admits an identifying code, we
say that $G$ is identifiable and denote by $\gamma^{ID}(G)$
 the minimum cardinality of an
identifying code of $G$. In this paper, we study the identifying
code of the corona product $H\odot G$ of graphs $H$ and $G$. We
first give a necessary and sufficient condition  for the
identifiable corona product $H\odot G$, and then express
$\gamma^{ID}(H\odot G)$ in terms of $\gamma^{ID}(G)$ and the (total)
domination number of $H$. Finally, we compute $\gamma^{ID}(H\odot
G)$ for some special graphs $G$.

\medskip



\end{abstract}
\bigskip

\bigskip

\section{Introduction}

Let $G$ be a finite graph. We often denote by $V(G)$ the vertex set
of $G$.
For $x\in V(G)$,
the {\em neighborhood} $N_G(x)$ of $x$ is the set of vertices
adjacent to $x$; the {\em closed neighborhood} $N_G[x]$ of $x$ is
the union of $\{x\}$ and $N_G(x)$. For subsets $C$ and $S$ of
$V(G)$, we say that $C$ {\em covers} $S$ if the set $N_G[x]\cap C$
is nonempty for each $x\in S$; we say that $C$ {\em separates} $S$
if the sets $N_G[x]\cap C$ are distinct for all $x\in S$. An {\em
identifying code}  of $G$ is a set of vertices which covers and
separates $V(G)$. If $G$ admits an identifying code, we say that $G$
is \emph{identifiable} and denote by $\gamma^{ID}(G)$ the minimum
cardinality of an identifying code of $G$. Note that $G$ is
identifiable if and only if the sets $N_G[x]$ are distinct for all
$x\in V(G)$.

The concept of identifying codes was  introduced by Karpovsky et al.
\cite{Ka1} to model a fault-detection problem in multiprocessor
systems.  It was noted in \cite{Ch,Coh} that determining the
identifying code with the minimum cardinality in a graph is an
NP-complete problem. Many researchers  focused on studying
identifying codes of some restricted graphs, for example, paths
\cite{Be}, cycles \cite{Be,Gr,Xu}, grids \cite{ben,coh2,ho1} and
triangle-free graphs \cite{fo}. The identifying codes of graph
products were studied; see  \cite{Bl1,Gra,Ho,ja,Ka2,Mo} for
cartesian products,  \cite{fe} for lexicographic products  and
\cite{ra} for  direct products.

The {\em corona product} $H\odot G$ of two graphs $H$ and $G$ is
defined as the graph obtained from $H$ and $G$ by taking one copy of
$H$ and $|V(H)|$ copies of $G$ and joining by an edge each vertex
from the $i$th-copy of $G$ with the $i$th-vertex of $H$. For each
$v\in V(H)$, we often refer to $G_v$ the copy of $G$ connected to
$v$ in $H\odot G$.

This paper is aimed to investigate identifying codes of the corona
product $H\odot G$ of graphs $H$ and $G$. In Section 2, we first
give a necessary and sufficient condition  for the identifiable
corona product $H\odot G$, and then construct   some identifying
codes of $H\odot G$. In Section 3, some inequalities for
$\gamma^{ID}(H\odot G)$ are established. In Section 4, we express
$\gamma^{ID}(H\odot G)$ in terms of $\gamma^{ID}(G)$ and the (total)
domination number of $H$. In Section 5,  we compute
$\gamma^{ID}(H\odot G)$ for some special graphs $G$.

\section{Constructions}
In this section, we first give a necessary and sufficient condition
for the identifiable corona product $H\odot G$, and then construct
some identifying codes of $H\odot G$.

\begin{thm}\label{identifiable}
Let $G$ be a graph.

{\rm(i)} Suppose $K_1$ is a trivial graph. Then $K_1\odot G$ is
identifiable if and only if $G$ is an identifiable graph with
maximum degree  at most $|V(G)|-2$.

{\rm(ii)} If $H$ is a nontrivial connected  graph, then $H\odot G$
is identifiable if and only if $G$ is identifiable.
\end{thm}
\proof (i) Write $V(K_1)=\{v\}$. Note that
$N_{K_1\odot G}[v]=V(K_1\odot G)$.
For any vertices $x$ and $y$ of $G_v$, we have
$N_{K_1\odot G}[x]=N_{K_1\odot G}[y]$ if and only if $N_{G_v}[x]=N_{G_v}[y]$.
Hence, the desired result follows.

(ii) If $H\odot G$ is identifiable, then $G_v$ is identifiable for each $v\in V(H)$,
which implies that $G$ is identifiable.
Conversely, suppose that $G$ is identifiable.
Pick any two distinct vertices $x$ and $y$ of $H\odot G$.
If $\{x,y\}\not\subseteq V(G_v)$ for any $v\in V(H)$, then
$N_{H\odot G}[x]\neq N_{H\odot G}[y]$.
If there exists a vertex $v\in V(H)$ such that $\{x,y\}\subseteq V(G_v)$,
 by $N_{G_v}[x]\neq N_{G_v}[y]$ we have $N_{H\odot G}[x]\neq N_{H\odot G}[y]$.
So $H\odot G$ is identifiable.
$\qed$

In the remaining of this section, some identifying codes of
the identifiable corona product $H\odot G$ are constructed.
We begin by a useful lemma.

\begin{lemma}\label{corona-identifying}
A set $C$ of vertices in the corona product $H\odot G$ is an
identifying code if, for each $v\in V(H)$,  the following three
conditions hold.

{\rm(i)} $C\cap V(G_v)$ is nonempty and separates $V(G_v)$ in $G_v$.

{\rm(ii)} $N_H(v)\cap C\neq\emptyset$, or $C\cap V(G_v)\not\subseteq
N_{G_v}[x]$ for any $x\in V(G_v)$.

{\rm(iii)} $v\in C$, or $C\cap V(G_v)$ covers $V(G_v)$ in $G_v$.
\end{lemma}
\proof Since $C\cap V(G_v)\neq\emptyset$, the set $C\cap V(G_v)$ covers $\{v\}$.
Since $\{v\}$ covers $V(G_v)$, by (iii) the set $C\cap(V(G_v)\cup\{v\})$ covers $V(G_v)$.
It follows that $C$ covers $V(H\odot G)$.
Hence, we only need to show that, for any two distinct vertices $x$ and $y$
in $V(H\odot G)$,
\begin{equation}\label{c1}
N_{H\odot G}[x]\cap C\neq N_{H\odot G}[y]\cap C.
\end{equation}

{\em Case 1.} $\{x,y\}\cap V(H)\neq\emptyset$. Without loss of
generality, assume that $x\in V(H)$. If $y\in V(H\odot G)\setminus
V(G_x)$, pick $z\in C\cap V(G_x)$, then $z\in (N_{H\odot G}[x]\cap
C)\setminus N_{H\odot G}[y]$, which implies that (\ref{c1}) holds.
Now suppose that $y\in V(G_x)$. If $C\cap V(G_x)\not\subseteq
N_{G_x}[y]$, then $N_{H\odot G}[x]\cap C\not\subseteq N_{H\odot
G}[y]$, and so (\ref{c1}) holds; If $C\cap V(G_x)\subseteq
N_{G_x}[y]$, by (ii) we can pick $z'\in N_H(x)\cap C$. Then $z'\in
(N_{H\odot G}[x]\cap C)\setminus N_{H\odot G}[y]$, and so (\ref{c1})
holds.

{\em Case 2.} $\{x,y\}\cap V(H)=\emptyset$. Then there exist
vertices $u$ and $v$ of $H$ such that $x\in V(G_u)$ and $y\in
V(G_v)$. If $u=v$, since $C\cap V(G_u)$ separates $\{x,y\}$ in
$G_u$, the set $C$ separates $\{x,y\}$ in $H\odot G$, which implies
that (\ref{c1}) holds; If $u\neq v$, then $N_{H\odot G}[x]\cap
N_{H\odot G}[y]=\emptyset$. Since $C$ covers $\{x,y\}$, the
inequality (\ref{c1}) holds. $\qed$


Next we shall construct identifying codes of $H\odot G$.

\begin{cor}\label{cons1}
Let $H$ be an arbitrary  graph and $G$ be an identifiable graph with
maximum degree at most $V(G)-2$. Then
$$
\bigcup_{v\in V(H)}S_v
$$
is an identifying code of $H\odot G$, where $S_v$ is an identifying
code of $G_v$ such that $S_v\not\subseteq N_{G_v}[x]$ for any vertex
$x$ of $G_v$.
\end{cor}
\proof It is immediate from Lemma~\ref{corona-identifying}.$\qed$

\begin{prop}\label{lemma1}
Let $S$ be a set of vertices in an identifiable graph $G$. If $S$ separates $V(G)$,
 then there exists a vertex $z\in V(G)$ such that
$S\cup\{z\}$ is an identifying code of $G$,
 and so $|S|\geq\gamma^{ID}(G)-1$.
\end{prop}
\proof If $S$ covers $V(G)$, then $S\cup\{z\}$ is an identifying
code of $G$ for any $z\in V(G)$. Now suppose that $S$ does not cover
$V(G)$. Then there exists a unique vertex $z\in V(G)$ such that
$N_G[z]\cap S=\emptyset$, which implies that $S\cup\{z\}$ is an
identifying code of $G$, as desired. $\qed$

From the above proposition, a set of vertices that separates the
vertex set is an identifying code, or is obtained from an
identifying code by deleting a  vertex. Now we use this set  of
vertices in $G$ and the vertex set of $H$ to construct identifying
codes of $H\odot G$.

\begin{cor}\label{cons2}
Let $H$ be a nontrivial connected graph  and $G$ be a nontrivial
identifiable graph.  Write
$$
C=\bigcup_{v\in V(H)}S_v\cup V(H),
$$
where  $S_v$ is a set of vertices separating $V(G_v)$ in $G_v$. Then
$C$ is an identifying code of $H\odot G$.
\end{cor}
\proof For each $v\in V(H)$, we have $C\cap V(G_v)=S_v\neq\emptyset$, $N_H(v)\cap C\neq\emptyset$ and $v\in C$.
It follows from Lemma~\ref{corona-identifying} that $C$ is an identifying code
of $H\odot G$.
$\qed$

Let $H$ be a graph. For a set $D$ of vertices, we say that $D$ is a
{\em dominating set} of $H$ if $D$ covers $V(H)$; we say that $D$ is
a {\em total dominating set} of $H$ if the set $N_H(x)\cap D$ is
nonempty for each $x\in V(H)$. The {\em domination number} of $H$,
denoted by $\gamma(H)$, is the minimum cardinality of a dominating
set of $H$; the {\em total domination number} of $H$, denoted by
$\gamma_t(H)$, is the minimum cardinality of a total dominating set
of $H$. Domination and its variations in graphs are now well
studied. The literature on this subject has been surveyed and
detailed in the the book \cite{ha1}.

The (total) dominating set of $H$ can be used to construct identifying codes of $H\odot G$.
The proofs of the following corollaries are immediate from Lemma~\ref{corona-identifying}.

\begin{cor}\label{cons3}
Let $H$ be an arbitrary graph and $G$ be an identifiable graph with
maximum degree at most $|V(G)|-2$. Suppose that $D$ is a dominating
set of $H$. Then
$$
\bigcup_{v\in V(H)}S_v\cup D
$$
is an identifying code of $H\odot G$, where  $S_v$ is an identifying
code of $G_v$ if $v\in V(H)\setminus D$;  $S_v$ is a set of vertices
separating $V(G_v)$ in $G_v$ such that $S_v\not\subseteq N_{G_v}[x]$
for any vertex $x$ of $G_v$ if $v\in D$.
\end{cor}

\begin{cor}\label{cons4}
Let $H$ be a nontrivial connected graph and $G$ be an identifiable
graph. Suppose that $T$ be a total dominating set of $H$. Then
$$
\bigcup_{v\in V(H)}S_v\cup T
$$
is an identifying code of $H\odot G$, where $S_v$ is an identifying
code of $G_v$.
\end{cor}

\section{Upper and lower bounds}
In this section, we shall establish some inequalities for
$\gamma^{ID}(H\odot G)$ by discussing the existence of some special
identifying codes of $G$.

In order to obtain upper bounds for $\gamma^{ID}(H\odot G)$, it
suffices to construct identifying codes of $H\odot G$.  By
Corollaries~\ref{cons1}, \ref{cons2} and \ref{cons3}, we need to
consider the identifying codes $S$ of $G$ satisfying one of the
following conditions:
\begin{itemize}
  \item [(a)] $|S|=\gamma^{ID}(G)$ and $S\not\subseteq N_G[x]$ for any $x\in
V(G)$.
  \item [(b)] $|S|=\gamma^{ID}(G)$ and there is a vertex $z\in S$ such that $S\setminus\{z\}$
separates $V(G)$.
  \item [(c)] $|S|=\gamma^{ID}(G)+1$ and there exists a vertex $z\in S$ such that
$S\setminus\{z\}$ separates $V(G)$ and
$S\setminus\{z\}\not\subseteq N_G[x]$ for any $x\in V(G)$.
\end{itemize}

The identifying codes satisfying (b) or (c) were studied in
\cite{Bl1,fe}.

\begin{lemma}\label{lemma4}
Let $G$ and $H$ be two graphs.
 If there exists an identifying code $S$ of $G$ satisfying {\rm(a)}, then
$
\gamma^{ID}(H\odot G)\leq|V(H)|\cdot\gamma^{ID}(G).
$
\end{lemma}
\proof For each $v\in V(H)$,
 let $S_v$ be the copy of $S$ in $G_v$.  Corollary~\ref{cons1} implies that
  $\cup_{v\in V(H)}S_v$ is an identifying code of
 $H\odot G$ with size $|V(H)|\cdot\gamma^{ID}(G)$, as desired.
$\qed$

\begin{lemma}\label{lemma6}
Let $G$ and $H$ be  two nontrivial graphs. Suppose that
 $H$ is connected.
If there is an identifying code $S$ of $G$ satisfying {\rm(b)}, then
$ \gamma^{ID}(H\odot G)\leq |V(H)|\cdot\gamma^{ID}(G). $
\end{lemma}
\proof 
Note that there exists a vertex $z\in S$ such that $S\setminus\{z\}$
separates $V(G)$. For each $v\in V(H)$, let $S_v$ be the copy of
$S\setminus\{z\}$ in $G_v$. It follows from Corollary~\ref{cons2}
that $\cup_{v\in V(H)}S_v\cup V(H)$ is an identifying code of
$H\odot G$ with size $|V(H)|\cdot\gamma^{ID}(G)$. Therefore, the
desired inequality holds. $\qed$

\begin{lemma}\label{lemma10}
Let $G$ and $H$ be two nontrivial graphs. If there exists an
identifying code $S$ of $G$ satisfying {\rm(c)}, then $
\gamma^{ID}(H\odot G)\leq|V(H)|\cdot\gamma^{ID}(G)+\gamma(H). $
\end{lemma}
\proof  Observe  that there exists a vertex $z\in S$ such that
$S\setminus\{z\}$ separates $V(G)$ and $S\setminus\{z\}\not\subseteq
N_G[x]$ for any vertex $x\in V(G)$. Suppose that $W$ is an
identifying code of $G$ with size $\gamma^{ID}(G)$ and $D$ is a
dominating set of $H$ with size $\gamma(H)$. For each $v\in D$, let
$S_v$ be the copy of $S\setminus\{z\}$ in $G_v$. For each $v\in
V(H)\setminus D$, let $S_v$ be the copy of $W$ in $G_v$. It follows
from Corollary~\ref{cons3} that $\cup_{v\in V(H)}S_v\cup D$ is an
identifying code of
 $H\odot G$ with size $|V(H)|\cdot\gamma^{ID}(G)+\gamma(H)$,
as desired.
$\qed$

With reference to Corollary~\ref{cons4}, let $T$ and $S_v$  have the sizes $\gamma_t(H)$ and
$\gamma^{ID}(G)$,
respectively. Then we get the following result immediately.

\begin{lemma}\label{lemma12}
Let  $G$ be an identifiable graph and $H$ be a nontrivial connected
graph. Then $ \gamma^{ID}(H\odot
G)\leq|V(H)|\cdot\gamma^{ID}(G)+\gamma_t(H). $
\end{lemma}

In the remaining of this section, we give lower bounds for
$\gamma^{ID}(H\odot G)$. We begin by discussing  the properties of
an identifying code of $H\odot G$.

\begin{lemma}\label{lemma2}
Let $C$ be an identifying code of $H\odot G$ and let $v$ be a vertex of the first factor $H$. Then $C\cap V(G_v)$ separates $V(G_v)$ in $G_v$. Moreover, if $v\not\in C$, then $C\cap V(G_v)$ is an identifying code of $G_v$.
\end{lemma}
\proof Note that $v$ is adjacent to every vertex in $V(G_v)$, and
there are no edges joining $V(H\odot G)\setminus(\{v\}\cup G_v)$
with $V(G_v)$. Since $C$ separates $V(G_v)$ in $H\odot G$, the set
$C\cap V(G_v)$ separates $V(G_v)$ in $G_v$. If $v\not\in C$, since
$C$ covers $V(G_v)$ in $H\odot G$, the set $C\cap V(G_v)$ covers
$V(G_v)$ in $G_v$, which implies that $C\cap V(G_v)$ is an
identifying code of $G_v$. $\qed$

\begin{lemma}\label{lemma3}
If $H\odot G$ is identifiable, then
$
\gamma^{ID}(H\odot G)\geq |V(H)|\cdot\gamma^{ID}(G).
$
\end{lemma}
\proof Let $C$ be an identifying code of $H\odot G$ with size $\gamma^{ID}(H\odot G)$.
Combining Lemma~\ref{lemma2} and Proposition~\ref{lemma1},
 we have
\begin{equation*}
|C\cap V(G_v)|\geq\left\{
\begin{array}{ll}
\gamma^{ID}(G)-1,&\textup{if } v\in V(H)\cap C,\\
\gamma^{ID}(G),&\textup{if }v\in V(H)\setminus C.
\end{array}\right.
\end{equation*}
Then
$$
\gamma^{ID}(H\odot G)=\sum_{v\in V(H)\cap C}(|C\cap V(G_v)|+1)+\sum_{v\in V(H)\setminus C}|C\cap V(G_v)|\geq |V(H)|\cdot\gamma^{ID}(G),
$$
as desired.
$\qed$

\begin{lemma}\label{lemma5}
Let $G$ be an identifiable graph with maximum degree at most $|V(G)|-2$. If
any identifying code of $G$ does not satisfy {\rm(a)},
then
$
\gamma^{ID}(K_1\odot G)\geq\gamma^{ID}(G)+1.
$
\end{lemma}
\proof By Theorem~\ref{identifiable}, the coronal product $K_1\odot G$ is identifiable.
Hence, Lemma~\ref{lemma3} implies that $\gamma^{ID}(K_1\odot G)\geq\gamma^{ID}(G)$. Suppose for
the contradiction that there exists an identifying code $C$ of $K_1\odot G$ with size $\gamma^{ID}(G)$.
Write $V(K_1)=\{v\}$.

{\em Case 1.} $v\not\in C$. Then $C$ is an identifying code of $G_v$
with cardinality $\gamma^{ID}(G)$ by Lemma~\ref{lemma2}. Hence,
there is a vertex $x\in V(G_v)$ such that $C\subseteq N_{G_v}[x]$,
which implies that $N_{H\odot G}[x]\cap C=C=N_{H\odot G}[v]\cap C$,
a contradiction.

{\em Case 2.} $v\in C$. Then $C\cap V(G_v)=C\setminus\{v\}$.
Combining Proposition~\ref{lemma1} and Lemma~\ref{lemma2}, there
exists a vertex $z\in V(G_v)$ such that $(C\setminus\{v\})\cup\{z\}$
is an identifying code of $G_v$ with cardinality $\gamma^{ID}(G)$.
Hence, we have $(C\setminus\{v\})\cup\{z\}\subseteq N_{G_v}[y]$ for
some $y\in V(G_v)$, which implies that $(C\setminus\{v\})\subseteq
N_{G_v}[y]$. Consequently, we get $N_{H\odot G}[y]\cap C=C=N_{H\odot
G}[v]\cap C$, a contradiction. $\qed$

\begin{lemma}\label{lemma7}
Suppose that $C$ is an identifying code of $H\odot G$. If
any identifying code of $G$ does not satisfy {\rm(b)},
then $|C\cap V(G_v)|\geq\gamma^{ID}(G)$ for each $v\in V(H)$.
\end{lemma}
\proof  Lemma~\ref{lemma2} implies that $C\cap V(G_v)$ separates $V(G_v)$ in $G_v$. Then
$|C\cap V(G_v)|\geq\gamma^{ID}(G)-1$ by Proposition~\ref{lemma1}.
If $|C\cap V(G_v)|=\gamma^{ID}(G)-1$,  there exists
a vertex $z\in V(G)$ such that $(C\cap V(G_v))\cup\{z\}$ is an identifying code
of $G_v$ satisfying (b), a contradiction.
$\qed$

For a set $C$ of vertices in $H\odot G$, write
$$
H(C)=V(H)\cap C,\quad H'(C)=\{v\mid v\in V(H), |C\cap G_v|\geq\gamma^{ID}(G)+1\}.
$$

\begin{lemma}\label{lemma8}
Suppose that $C$ is an identifying code of  $H\odot G$.
If any identifying code  of $G$ does not satisfy {\rm(b)}, then
$
|C|\geq|V(H)|\cdot\gamma^{ID}(G)+|H(C)|+|H'(C)|.
$
\end{lemma}
\proof Write $H_1=V(H)\setminus(H(C)\cup H'(C))$, $H_2=H'(C)\setminus H(C)$,
$H_3=H(C)\setminus H'(C)$ and $H_4=H(C)\cap H'(C)$.
Let $C_v=C\cap V(G_v)$.
By Lemma~\ref{lemma7} we get $|C_v|=\gamma^{ID}(G)$ for
each $v\in H_1\cup H_3$.
Then
\begin{eqnarray*}
|C|&=&\sum_{v\in H_1}|C_v|+\sum_{v\in H_2}|C_v|+\sum_{v\in H_3}(|C_v|+1) +\sum_{v\in H_4}(|C_v|+1)\\
&\geq& |H_1|\gamma^{ID}(G)+|H_2|(\gamma^{ID}(G)+1)+|H_3|(\gamma^{ID}(G)+1)+|H_4|(\gamma^{ID}(G)+2)\\
&=&|V(H)|\cdot\gamma^{ID}(G)+|H(C)|+|H'(C)|,
\end{eqnarray*}
as desired.
$\qed$

\begin{lemma}\label{lemma9}
Let $G$ be a nontrivial identifiable graph and $H$ be a nontrivial
connected graph. If each identifying code of $G$ satisfies neither
{\rm(a)} nor {\rm(b)}, then $ \gamma^{ID}(H\odot
G)\geq|V(H)|\cdot\gamma^{ID}(G)+\gamma(H). $
\end{lemma}
\proof  Theorem~\ref{identifiable} implies that $H\odot G$ is identifiable.
Let $C$ be an identifying code of $H\odot G$ with size $\gamma^{ID}(H\odot G)$.
Write
$$
D=H(C)\cup H'(C).
$$
We shall show that $D$ is a dominating set of $H$. Pick any
$v\in V(H)\setminus D$. Note that $v\not\in C$ and $|C\cap V(G_v)|\leq\gamma^{ID}(G)$.
Then $C\cap V(G_v)$ is an identifying code of $G_v$ with size $\gamma^{ID}(G_v)$
by Lemma~\ref{lemma2}.
Since each identifying code of $G_v$ does not satisfy (a), there exists a vertex
$x\in V(G_v)$ such that $C\cap V(G_v)\subseteq N_{G_v}[x]$.
Since $N_{H\odot G}[v]\cap C\neq N_{H\odot G}[x]\cap C=C\cap V(G_v)$,
 we have $N_H(v)\cap H(C)\neq\emptyset$,
which implies that $N_H(v)\cap D\neq\emptyset$.
Then $D$ is a dominating set of $H$.

Hence, we have $|D|\geq\gamma(H)$.
By Lemma~\ref{lemma8}, we get
$$
\gamma^{ID}(H\odot G)=|C|\geq |V(H)|\cdot\gamma^{ID}(G)+|H(C)\cup H'(C)|\geq |V(H)|\cdot\gamma^{ID}(G)+\gamma(H),
$$
as desired.
$\qed$

\begin{lemma}\label{lemma11}
Let $G$ be a nontrivial identifiable graph and $H$ be a nontrivial
connected graph. If each identifying code of $G$ satisfies none of
the conditions {\rm(a)}, {\rm(b)} and {\rm(c)}, then $
\gamma^{ID}(H\odot G)\geq|V(H)|\cdot\gamma^{ID}(G)+\gamma_t(H). $
\end{lemma}
\proof For each vertex $v\in V(H)$, pick a vertex $v'\in N_H(v)$.
Theorem~\ref{identifiable} implies that $H\odot G$ is identifiable.
Let $C$ be an identifying code of $H\odot G$ with size
$\gamma^{ID}(H\odot G)$. Write
$$
T=H''(C)\cup H(C),
$$
where  $H''(C)=\{v'\mid v\in H'(C)\}$.

We claim that $T$ is a total dominating set of $H$. Pick any $v\in
V(H)$. If $v\in H'(C)$, since $N_H(v)\cap H''(C)\neq\emptyset$ we
have $N_H(v)\cap T\neq\emptyset$. Now suppose that $v\not\in H'(C)$.
By Lemma~\ref{lemma7} we get $|C\cap V(G_v)|=\gamma^{ID}(G_v)$. If
$C\cap V(G_v)\not\subseteq N_{G_v}[x]$ for any vertex $x\in V(G_v)$,
then $C\cap V(G_v)$ is not an identifying code of $G_v$. It follows
from Lemma~\ref{lemma2} and Proposition~\ref{lemma1} that there
exists a vertex $z\in V(G_v)$ such that $(C\cap V(G_v))\cup\{z\}$ is
an identifying code of $G_v$ satisfying (c),  a contradiction.
Therefore, there exists a vertex $x\in V(G_v)$ such that $C\cap
V(G_v)\subseteq N_{G_v}[x]$. Since $N_{H\odot G}[v]\cap C\neq
N_{H\odot G}[x]\cap C$, we have $N_H(v)\cap H(C)\neq\emptyset$,
which implies that $N_H(v)\cap T\neq\emptyset$. Hence, our claim is
valid.

Since $|T|\geq\gamma_t(H)$ and $|H'(C)|\geq |H''(C)|$, we get $|H'(C)|+|H(C)|\geq\gamma_t(H)$.
By Lemma~\ref{lemma8}, we have
$$
\gamma^{ID}(H\odot G)=|C|\geq  |V(H)|\cdot\gamma^{ID}(G)+\gamma_t(H),
$$
as desired.
$\qed$

\section{Minimum cardinality}
In this section, we shall compute $\gamma^{ID}(H\odot G)$.


\begin{thm}\label{main3}
Let $G$ and $H$ be two nontrivial  graphs. Suppose that $H$ is
connected. If there exists an identifying code of $G$ satisfying
{\rm(a)} or {\rm(b)}, then
$$
\gamma^{ID}(H\odot G)=|V(H)|\cdot\gamma^{ID}(G).
$$
\end{thm}
\proof
It is immediate from Theorem~\ref{identifiable}, Lemmas~\ref{lemma4}, \ref{lemma6} and \ref{lemma3}.
$\qed$

\begin{thm}\label{main4}
Let $G$ be a nontrivial identifiable graph and $H$ be a nontrivial
connected graph. Suppose that each identifying code of $G$ satisfies
neither {\rm(a)} nor {\rm(b)}.

{\rm(i)} If there exists an identifying code of $G$ satisfying {\rm(c)}, then
$$
\gamma^{ID}(H\odot G)=|V(H)|\cdot\gamma^{ID}(G)+\gamma(H).
$$

{\rm(ii)} If any identifying code  of $G$ does not satisfy {\rm(c)}, then
$$
\gamma^{ID}(H\odot G)=|V(H)|\cdot\gamma^{ID}(G)+\gamma_t(H).
$$
\end{thm}
\proof
 (i) holds by Lemmas~\ref{lemma10} and \ref{lemma9}.
By Lemmas~\ref{lemma12} and \ref{lemma11}, (ii) holds. $\qed$

Now, we compute $\gamma^{ID}(K_1\odot G)$ and $\gamma^{ID}(H\odot
K_1)$.

\begin{thm}\label{main1}
Suppose that $G$ is an identifiable graph with maximum
degree  at most $|V(G)|-2$.

{\rm (i)} If there exists an identifying code of $G$ satisfying {\rm(a)}, then
$$
\gamma^{ID}(K_1\odot G)=\gamma^{ID}(G).
$$

{\rm(ii)} If any identifying code of $G$ does not satisfy {\rm(a)}, then
$$
\gamma^{ID}(K_1\odot G)=\gamma^{ID}(G)+1.
$$
\end{thm}
\proof Theorem~\ref{identifiable} implies that $K_1\odot G$ is identifiable.

(i) It is immediate from Lemmas~\ref{lemma4} and \ref{lemma3}.

(ii) By Lemma~\ref{lemma5} we only need to construct an identifying code of $K_1\odot G$
with size $\gamma^{ID}(G)+1$.
Let $W$ be an identifying code of $G$ with size $\gamma^{ID}(G)$.
Note that there exists a unique vertex $x\in V(G)$ such that
$W\subseteq N_G[x]$. Pick $y\in V(G)\setminus N_G[x]$. Write $V(K_1)=\{v\}$.
Let $S_v$ be the copy of $W\cup\{y\}$ in $G_v$.
Then $S_v$ is an identifying code of $G_v$
with $S_v\not\subseteq N_{G_v}[z]$ for any vertex $z\in V(G_v)$.
It follows from Corollary~\ref{cons1} that $S_v$ is an identifying code of $K_1\odot G$
with size $\gamma^{ID}(G)+1$, as desired.
$\qed$

\begin{cor}\label{cor}
Let $G$ be an identifiable graph and $H$ be a connected graph.
Suppose that $G$ satisfies one of the following conditions.

{\rm(i)} The graph $G$ is not connected.

{\rm(ii)} The diameter of $G$ is at least five.

{\rm(iii)} The maximum degree of $G$ is less than
$\gamma^{ID}(G)-1$.
\\*
Then
$$
\gamma^{ID}(H\odot G)=|V(H)|\cdot\gamma^{ID}(G).
$$
\end{cor}
\proof Note that the identifying codes of $G$ with size
$\gamma^{ID}(G)$ satisfy (a). Combining Theorems~\ref{main3} and
\ref{main1}, we get the desired result. $\qed$

\begin{thm}\label{kn}
Let $n\geq 2$. Then $\gamma^{ID}(K_n\odot K_1)=n+1$, where  $K_n$ is
the complete graph on $n$ vertices.
\end{thm}
\proof
Since $K_1$ is identifiable, Theorem~\ref{identifiable} implies that
$K_n\odot K_1$ is identifiable.
Write $V=V(K_n)=\{v_1,\ldots,v_n\}$.
For each $i\in\{1,\ldots,n\}$, denote by $\{u_i\}$ the vertex set of the copy of $K_1$
connected to $v_i$ in $K_n\odot K_1$. Write $V'=\{u_1,\ldots,u_n\}$.
Note that $V(K_n\odot K_1)=V\cup V'$.
Let $C$ be an identifying code of $K_n\odot K_1$ with size $\gamma^{ID}(K_n\odot K_1)$.
We have the following two claims.

{\em Claim 1.} $|V\cap C|\geq 2$. In fact, for any $i\in\{1,\ldots,n\}$, since
$$
(V\cup\{u_i\})\cap C=N_{K_n\odot K_1}[v_i]\cap C\neq N_{K_n\odot K_1}[u_i]\cap C=\{u_i,v_i\}\cap C,
$$
 we have
$|(V\setminus\{v_{i}\})\cap C|\geq 1$. So $|V\cap C|\geq 2$.

{\em Claim 2.}  $|V'\cap C|\geq n-1$. In fact, if there exist two
distinct vertices $u_{i}$ and $u_{j}$ neither of  which   belongs to
$C$, then $N_{K_n\odot K_1}[v_{i}]\cap C=N_{K_n\odot K_1}[v_{j}]\cap
C$, a contradiction.

Combining Claim 1 and Claim 2, we have
$$
\gamma^{ID}(K_n\odot K_1)=|V\cap C|+|V'\cap C|\geq n+1.
$$
It is routine to show that $ \{u_{i}\mid 2\leq i\leq
n\}\cup\{v_{1},v_{2}\} $ is an identifying code of $K_n\odot K_1$
with size $n+1$. Hence, the desired result follows. $\qed$

\begin{thm}\label{main2}
Let $H$ be a connected graph that is not complete. Then
$$
\gamma^{ID}(H\odot K_1)=|V(H)|.
$$
\end{thm}
\proof Theorem~\ref{identifiable} implies that $H\odot K_1$ is
identifiable. Since $\gamma^{ID}(K_1)=1$, by Lemma~\ref{lemma3} it
suffices to construct an identifying code of $H\odot K_1$ with size
$|V(H)|$.

 For any $u, v\in V(H)$, define $u\equiv v$ if
$N_{H}(u)=N_{H}(v)$. Note that $``\equiv"$ is an equivalence
relation. Let $O_u$ denote the equivalence class containing $u$.
Pick a representative system $D$  with respect to this equivalence
relation. For each $v\in V(H)$, denote by $\{u_v\}$ the vertex set
of the copy of $K_1$ connected to $v$ in $H\odot K_1$. Let
$$
C=\{u_{v}\mid v\in V(H)\setminus D\}\cup D.
$$
Observe that $|C|=|V(H)|$. Since $C$ covers $V(H\odot K_1)$, it
suffices to show that, for any two distinct vertices $x$ and $y$ of
$H\odot K_1$,
\begin{equation}\label{c2}
N_{H\odot K_1}[x]\cap C\neq N_{H\odot K_1}[y]\cap C.
\end{equation}

{\em Case 1.} $|\{x,y\}\cap V(H)|=2$. If $N_H[x]\neq N_H[y]$, there
exists a vertex $z\in V(H)$ such that $\{z\}$ separates $\{x,y\}$ in
$H$. Note that there exists a vertex $z'\in D$ such that
$O_{z'}=O_z$. Then $N_H[z']=N_H[z]$, and so $\{z'\}$ separates
$\{x,y\}$ in $H$. It follows that $\{z'\}$ separates $\{x,y\}$ in
$H\odot K_1$. Since $z'\in C$, the inequality (\ref{c2}) holds. If
$N_H[x]=N_H[y]$, then $O_x=O_y$, which implies that $x\not\in D$ or
$y\not\in D$. Without loss of generality, we may assume that
$x\not\in D$. Then $u_{x}\in (N_{H\odot K_1}[x]\cap C)\setminus
N_{H\odot K_1}[y]$, which implies that (\ref{c2}) holds.

{\em Case 2.} $|\{x,y\}\cap V(H)|=1$. Without loss of generality, assume that
$x\in V(H)$ and $y=u_v$ for some $v\in V(H)$.
If $x\neq v$, since both $\{x\}$ and $\{u_{x}\}$ separate $\{x,y\}$ in $H\odot G$,
we obtain (\ref{c2}) by $\{x,u_{x}\}\cap C\neq\emptyset$.
Now suppose that $x=v$. Since $H$ is not complete, we have $|D|\geq 2$.
Hence, there is a vertex $w\in D$ such that $w$ is adjacent to $x$ in $H$.
It follows that $w\in (N_{H\odot K_1}[x]\cap C)\setminus N_{H\odot K_1}[y]$,
and so (\ref{c2}) holds.

{\em Case 3.} $|\{x,y\}\cap V(H)|=0$. Then $N_{H\odot K_1}[x]\cap
N_{H\odot K_1}[y]=\emptyset$. Since $C$ covers $\{x,y\}$ in $H\odot
G$, the inequality (\ref{c2}) holds. $\qed$

Let $T_1=K_1\odot K_1$ and $T_n=T_{n-1}\odot K_1$ for $n\geq 2$. We
call $T_n$ a {\em binomial tree}, which is a useful data structure
in the context of algorithm analysis and design \cite{co}. Note that
$T_n$ is a spanning tree of the hypercube $Q_n$. The problem of
computing $\gamma^{ID}(Q_n)$ is still open. By Theorem~\ref{main2},
we get the following corollary.

\begin{cor}
Let $n\geq 3$. Then $\gamma^{ID}(T_n)=2^{n-1}$.
\end{cor}

 For a connected
graph with pendant edges, we have the following more general result
than Theorem~\ref{main2}.

\begin{cor}
Let $H$ be a connected graph with $m$ vertices. Suppose that $H_1$
is a graph obtained from $H$ by adding $n_i (\geq 1)$ pendant edges
to the $i$th-vertex of $H$. If $H_1$ is not isomorphic to $K_m\odot
K_1$, then
\begin{equation}\label{c4}
\gamma^{ID}(H_1)=\sum_{i=1}^mn_i.
\end{equation}
\end{cor}
\proof It is routine to show that (\ref{c4}) holds for $m=1$. Now
suppose $m\geq 2$. Write $V(H)=\{v_1,\ldots,v_m\}$. For each
$i\in\{1,\ldots,m\}$, let $S_i=\{u_{ij}\mid 1\leq j\leq n_i\}$ be
the set of vertices adjacent to $v_i$ in $V(H_1)\setminus V(H)$.
Then the subgraph of $H_1$ induced by $S_i$ is isomorphic to
$\overline K_{n_i}$. Similar to the proof of Lemma~\ref{lemma3}, we
have
\begin{equation*}
\gamma^{ID}(H_1)\geq \sum_{i=1}^m\gamma^{ID}(\overline
K_{n_i})=\sum_{i=1}^m n_i.
\end{equation*}
In order to prove (\ref{c4}), it suffices to construct an
identifying code of $H_1$ with size $\sum_{i=1}^m n_i$.

{\em Case 1.} $H$ is a complete graph. Then there exists an index
$j\in\{1,\ldots,m\}$ such that $n_j\geq 2$. Pick
$k\in\{1,\ldots,m\}\setminus\{j\}$. It is routine to show that
$$
\{v_j,v_k\}\cup(S_j\setminus\{u_{j1}\})\cup(S_k\setminus\{u_{k1}\})\cup
\bigcup_{i\in\{1,\ldots,n\} \setminus \{j,k\}} S_i
$$
is an identifying code of $H_1$ with size $\sum_{i=1}^m n_i$.

 {\em Case 2.} $H$ is not a complete graph. Write
$$
A=\bigcup_{i=1}^m\{v_i,u_{i1}\},\quad B=V(H_1)\setminus A.
$$
Then the subgraph $H_1[A]$ of $H_1$ induced by $A$ is isomorphic to
$H\odot K_1$. Pick a subset $A_0\subseteq A$ such that  $A_0$ is an
identifying code of $H_1[A]$ with the minimum cardinality. By
Theorem~\ref{main2} we have $|A_0|=\gamma^{ID}(H\odot K_1)=m$. Let $
C=A_0\cup B. $ Note that $|C|=\sum_{i=1}^mn_i$. It suffices to show
that $C$ is an identifying code of $H_1$. The fact that $A_0$ covers
$A$ in $H_1[A]$ implies that $C$ covers $V(H_1)$ in $H_1$.
Therefore, we only need to show that, for any two distinct vertices
$x$ and $y$ of $H_1$,
\begin{equation}\label{c3}
N_{H_1}[x]\cap C\neq N_{H_1}[y]\cap C.
\end{equation}

{\em Case 2.1.} $\{x,y\}\subseteq A$. Then there is a vertex $z\in
A_0$ such that $\{z\}$ separates $\{x,y\}$ in $H_1[A]$, which
implies that $z\in C$ and $\{z\}$ separates $\{x,y\}$ in $H_1$. So
(\ref{c3}) holds.

{\em Case 2.2.} $\{x,y\}\not\subseteq A$. Without loss of
generality, we may assume that $x\not\in A$. Then $x\in B$. Write
$x=u_{ij}$, where $1\leq i\leq m$ and $2\leq j\leq n_i$. If $y\neq
v_i$, then $x\in (N_{H_1}[x]\cap C)\setminus N_{H_1}[y]$, which
implies that (\ref{c3}) holds. Now suppose that $y=v_i$. Since
$\{u_{i1},v_i\}\subseteq A$, there exists a vertex $z\in A_0$ such
that $\{z\}$ separates $\{u_{i1},v_i\}$ in $H_1[A]$, which implies
that $z\in C$ and $\{z\}$ separates $\{x,y\}$ in $H_1$. So
(\ref{c3}) holds. $\qed$

\section{Examples}
In  this section, we shall find some graphs satisfying each
condition  in Theorems~\ref{main3}, \ref{main4} and \ref{main1},
respectively. As a result, we compute $\gamma^{ID}(H\odot G)$ for
some special graphs $G$.

The minimum cardinality of an identifying code of the path $P_n$ or the cycle $C_n$ was computed in \cite{Be, Gr}.

\begin{prop}\label{pc} {\rm (\cite{Be, Gr})} {\rm(i)} For $n\geq3$,
$\gamma^{ID}(P_{n})=\lfloor\frac{n}{2}\rfloor+1$.

{\rm(ii)} For $n\geq 6$, $\gamma^{ID}(C_{n})=\left\{
\begin{array}{ll}
\frac{n}{2}, &n ~\textup{is even},\\
\frac{n+3}{2}, &n ~\textup{is odd}.
\end{array}\right.$
\end{prop}

Note that $\gamma^{ID}(C_4)=\gamma^{ID}(C_5)=3$. Each identifying
code of $P_3$, $P_4$, $C_4$ or $C_5$ satisfies none of the
conditions (a), (b) and (c). There exists an identifying code  of
$P_n$ (resp. $C_n$) satisfying (a) for $n\geq 5 $ (resp. $n\geq 6$).
Combining Theorems~\ref{identifiable}, \ref{main3}, \ref{main4},
\ref{main1}, Corollary~\ref{cor} and Proposition~\ref{pc}, we get
Examples \ref{ex1}, \ref{ex2} and Corollary~\ref{pc2}.

\begin{example}\label{ex1}
Let $F_n$ be a fan, that is  $F_n=K_1\odot P_n$. If $1\leq n\leq 3$,
then $F_n$ is not identifiable;  If $n\geq 4$, then
$$
\gamma^{ID}(F_n)=\left\{
\begin{array}{ll}
4, &\textup{if } n =4,\\
\lfloor\frac{n}{2}\rfloor+1, &\textup{if }n \geq 5.
\end{array}\right.
$$
\end{example}

\begin{example}\label{ex2}
Let $W_n$ be a wheel, that is  $W_n=K_1\odot C_n$. Then $W_3$ is not
identifiable. For $n\geq 4$, we have
$$
\gamma^{ID}(W_n)=\left\{
\begin{array}{ll}
4& \textup{if } n=4,\\
\frac{n}{2}, &\textup{if }n \textup{ is even and }n\geq 6,\\
\frac{n+3}{2}, &\textup{if }n \textup{ is odd and }n\geq 5.
\end{array}\right.
$$
\end{example}

\begin{cor}\label{pc2}
Let $H$ be a nontrivial connected graph with $m$ vertices.

{\rm(i)} $\gamma^{ID}(H\odot P_3)=2m+\gamma_t(H).$

{\rm(ii)}
$\gamma^{ID}(H\odot P_4)=\gamma^{ID}(H\odot C_4)=\gamma^{ID}(H\odot C_5)=3m+\gamma_t(H).$

{\rm(iii)} For $n\geq 5$, we have
$\gamma^{ID}(H\odot P_n)=m(\lfloor\frac{n}{2}\rfloor+1)$.

{\rm(iv)} For $n\geq 6$, we have
$
\gamma^{ID}(H\odot C_n)=\left\{
\begin{array}{ll}
\frac{mn}{2}, &n ~\textup{is even},\\
\frac{m(n+3)}{2}, &n ~\textup{is odd}.
\end{array}\right.
$
\end{cor}

Let $S_n$ be a star, that is $S_n=K_1\odot\overline K_n$, where
$\overline K_n$ is the empty graph on $n$ vertices. Suppose $n\geq
3$. By Corollary~\ref{cor}, we get $\gamma^{ID}(S_n)=n$. Each
identifying code of $S_n$ with size $n$ satisfies (b). By
Theorem~\ref{main3}, we have the following result.

\begin{cor}\label{cor3}
Let $H$ be a nontrivial connected graph with $m$ vertices.
 If $n\geq 3$, then $\gamma^{ID}(H\odot S_n)=mn.$
\end{cor}

\begin{center}
\setlength{\unitlength}{1mm}
\begin{picture}(23,26)
\put(14,16){\line(0,1){8}} \put(14,16){\line(-2,-1){8}} \put(14,16){\line(2,-1){8}}
\put(14,8){\line(-2,1){8}} \put(14,8){\line(2,1){8}}
\put(6,12){\line(0,1){8}} \put(22,12){\line(0,1){8}}
\put(14,24){\line(-2,-1){8}} \put(14,24){\line(2,-1){8}}
\put(14,8){\circle*{1}}\put(13,4){1}
\put(6,12){\circle*{1}}\put(3,9){6}
\put(22,12){\circle*{1}}\put(23,9){2}
\put(14,16){\circle*{1}}\put(13,12){0}
\put(6,20){\circle*{1}}\put(3,20){5}
\put(22,20){\circle*{1}}\put(23,20){3}
\put(14,24){\circle*{1}}\put(13,25){4}
\put(3,0){Figure 1: $G_3$}
\end{picture}
\end{center}

Let $G_3$ be the graph in Figure 1. Note that $\gamma^{ID}(G_3)=3$
and each identifying code with size three is contained in $\{0, 2,
4, 6\}$. Any subset of $V(G_3)$ with size two can not separates
$V(G_3)$. Therefore, each identifying code of $G_3$ satisfies
neither (a) nor (b). The fact that $\{1,3,5\}$ separates $V(G_3)$
implies that $\{0,1,3,5\}$ is an identifying code of $G_3$
satisfying (c). By Theorems~\ref{main4}, we get the following
result.

\begin{cor}
Let $H$ be a nontrivial connected graph with $m$ vertices. Then
$$
\gamma^{ID}(H\odot G_3)=3m+\gamma(H).
$$
\end{cor}

\section*{Acknowledgement}
 This research is supported by NSFC(11271047) and the Fundamental Research
Funds for the Central University of China.

\end{CJK*}
\end{document}